\documentclass[12pt]{article}
\usepackage{graphicx}% Include figure files
\usepackage{amssymb}
\usepackage{amsfonts}
\usepackage{amsthm}
\usepackage{amsfonts}

\newcommand \disc {{\mathbb D}}

\newcommand{\complex}{{\mathbb C}}

\newcommand \ra {\rightarrow}

\newcommand{\ba}[1]{\begin{array}{#1}}
\newcommand{\ea}{\end{array}}

\newcommand{\be}{\begin{equation}}
\newcommand{\ee}{\end{equation}}
\newcommand{\bea}{\begin{eqnarray}}
\newcommand{\eea}{\end{eqnarray}}
\newcommand{\beann}{\begin{eqnarray*}}
\newcommand{\eeann}{\end{eqnarray*}}

\def\reff#1{(\ref{#1})}
\newtheorem{proposition}{Proposition}

\bibstyle{ams}

\begin{document}

\title{The Smart Kinetic Self-Avoiding Walk and \\ Schramm Loewner Evolution}

\author{Tom Kennedy
\\Department of Mathematics
\\University of Arizona
\\Tucson, AZ 85721
\\ email: tgk@math.arizona.edu
}

\maketitle 

\begin{abstract}
The smart kinetic self-avoiding walk (SKSAW) 
is a random walk which never intersects itself 
and grows forever when run in the full-plane.
At each time step the walk chooses 
the next step uniformly from among the allowable nearest neighbors 
of the current endpoint of the walk. In the full-plane a nearest neighbor is 
allowable if it has not been visited before and there is a path 
from the nearest neighbor to infinity through sites that have 
not been visited before. It is well known that on the hexagonal 
lattice the SKSAW in a bounded domain between two boundary points is 
equivalent to an interface in critical percolation, and hence its 
scaling limit is the chordal Schramm-Loewner evolution with $\kappa=6$
(SLE$_6$). Like SLE there are variants 
of the SKSAW depending on the domain and the initial and 
terminal points. 
On the hexagonal lattice these variants have been shown to converge
to the corresponding version of SLE$_6$.
It is believed that the scaling limit of all
these variants on any regular lattice is the corresponding 
version of SLE$_6$. We test this conjecture for the square lattice
by simulating the SKSAW in the full-plane 
and find excellent agreement with the predictions of full-plane SLE$_6$.
\end{abstract}

\newpage

\section{Introduction}

There are a variety of models that produce random walks that are 
simple, i.e., walks that do not have any self-intersections. 
In the {\it self-avoiding walk} (SAW) model, all simple walks with the 
same number of steps are given the same probability \cite{madras_slade}.
The scaling limit of this model is conjectured to be SLE$_{8/3}$
(Schramm-Loewner evolution with $\kappa=8/3$) 
\cite{lsw_saw}, and there is numerical evidence to support this conjecture
\cite{Kennedya,Kennedyb}.
{\it Loop-erased random walks} (LERW) are generated by erasing the loops in an
ordinary random walk in chronological order. The scaling limit of 
this model has been proved to be SLE$_2$ \cite{lerw_sle}.

In this paper we study another model of simple random walks in the plane
that was introduced in the physics literature in the mid 1980's under a 
couple of names. 
We will refer to it as the {\it smart kinetic self-avoiding walk} 
(SKSAW) for reasons we will explain later. 
The model can be defined on any lattice
in any number of dimensions, but our study will only be 
concerned with two dimensions. We are particularly interested in 
the relationship of the scaling limit of this model with  SLE$_6$. 

\begin{figure}[tbh]
\includegraphics{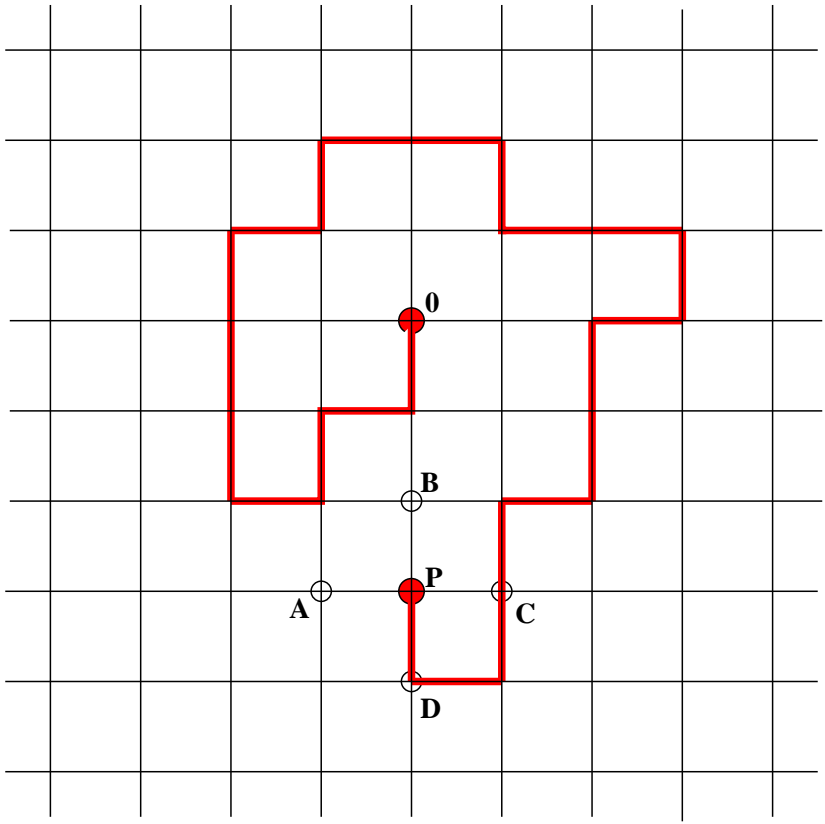}
\caption{\leftskip=25 pt \rightskip= 25 pt 
The walk starts at the origin $0$ and has been defined up to $P$. 
The nearest neighbors $C$ and $D$ are occupied. The nearest neighbor
$B$ is a trapping site. Only the nearest neighbor $A$ is allowed for 
the next step of the walk.
}
\label{fig_trap}
\end{figure}

There are several versions of the SKSAW corresponding to different choices
of the region in which the walk takes place and the choice of 
initial and terminal points.  
We first define one that we will call the ``full-plane'' SKSAW.
We start the walk at the origin, i.e., $\omega(0)=0$. 
For the first step $\omega(1)$ we choose one of the nearest neighbors of 
the origin with equal probability. 
Now suppose the steps $\omega(i)$ have been chosen for $i=0,1,\cdots,n$.
We say that a site $z$ is occupied if the walk has visited $z$ before,
i.e., $\omega(i)=z$ for some $0 \le i \le n$. We say that $z$ is 
a trapping site (with respect to $\infty$) if there does not exist 
a nearest neighbor path from $z$ to $\infty$ though 
unoccupied sites. Figure \ref{fig_trap} illustrates these definitions.
Now consider the nearest neighbors of $\omega(n)$
that are not occupied and are not trapping sites. 
We assign them equal probability and choose one at random to 
be $\omega(n+1)$. 
A simple induction argument shows there is always at least one 
such nearest neighbor. (Because $\omega(n)$ was not a trapping site
when it was chosen, it is connected to $\infty$.)
This process will run forever and generate a self-avoiding walk
in the plane from $0$ to $\infty$. So we will refer to it as the 
full-plane SKSAW from $0$ to $\infty$. 

The next version we consider is also defined in the full-plane and 
generates walks between two points. 
We take the origin to be one of the two points and denote the 
other point by $v$. For the moment we assume that $v$ is not the origin. 
We start the walk at the origin as before. Given the 
walk $\omega(i)$ for $i=0,1,2,\cdots,n$, we define occupied sites 
as before. We say that a site $z$ is a trapping site with respect to $v$
if there does not exist a nearest neighbor path from $z$ to $v$ though 
unoccupied sites. Now consider the nearest neighbors of $\omega(n)$
that are not occupied and are not trapping sites with respect to $v$. 
We choose one with equal probability. 
We stop the process when the walk reaches
$v$. It is possible that the walk never reaches $v$. We 
believe, but cannot prove, that this event has probability zero.
We refer to this model as the full-plane SKSAW from $0$ to $v$.

In the previous definition we can take $v$ to be the origin 
if we modify the definition slightly.
When the walk is at a nearest neighbor of the origin, we 
consider the origin to be unoccupied. So the origin
is an allowed choice for the next step. 
We stop the process when it returns to the origin. 
Again, there are walks that never return to the origin, but we 
expect this event to have probability zero. 
For the hexagonal lattice we will see later that this follows 
from a relation between the SKSAW and critical percolation.
Since this model will produce loops that pass through the origin, 
we will refer to it as the loop SKSAW.

The next version is the chordal SKSAW. We let $D$ be a simply 
connected domain and fix two lattice sites $u,v$ on its boundary.
(More precisely, we fix two sites on the boundary of $D$ and 
let $u,v$ be the closest lattice sites inside the domain.)  
We start the walk at $u$. We want the walk to end at $v$.
Given the walk $\omega(i)$ for $i=0,1,\cdots,n$ we define occupied sites
as before. We treat sites which are outside the domain the same 
as occupied sites - the walk is not allowed to visit them. 
A site $z$ is a trapping site (with respect to $v$) if there does not exit a
nearest neighbor path from $z$ to $v$ though sites which are inside the
domain and are not occupied.
Now consider the nearest neighbors of $\omega(n)$ which are  inside the 
domain, are not occupied and are not trapping sites. 
We pick one of these sites with equal probability.
(There will always be at least one such site.) We stop the process
when we reach $v$. If $D$ is bounded this will always happen. If 
$D$ is unbounded there can be walks that never reach $v$, but 
we believe this event has probability zero. 
If $D$ is unbounded we can take $v$ to be $\infty$. 
In this case the walk goes on forever and 
converges to $\infty$. 
We have assumed the domain is simply connected since we are primarily 
interested in the relationship of the scaling limit to SLE, but
the definition of the SKSAW does not require the 
domain to be simply connected.

The final version is the radial SKSAW. We now let $u$ be a site on the 
boundary of a simply connected domain $D$ 
and $v$ a site in the interior of the domain. The walk goes 
from $u$ to $v$ and is defined just as in the chordal case. 
Note that we can also run the walk from $v$ to $u$. 
By considering some simple examples, one finds that 
in general the distribution for the walk from $u$ to $v$ is not the 
same as that of the walk run from $v$ to $u$. An interesting question
is whether they have the same scaling limit. 

On a hexagonal lattice the chordal SKSAW in a domain $D$ 
is exactly equivalent to an interface in critical percolation on the 
hexagonal lattice in $D$ with suitable boundary conditions.
We review this well-known equivalence. Fix two sites $u,v$ on the 
boundary of $D$. We first color the hexagons around the boundary 
of $D$ in the following deterministic way. We start at $u$ 
and traverse the boundary in the clockwise direction and color the 
boundary hexagons black until we reach $v$. Then as we continue
in the clockwise direction around the boundary we color the hexagons
white until we reach $u$. The hexagons in the interior of the 
domain are randomly colored black or white with probability $1/2$. 
These boundary conditions force an interface which runs between $u$ and 
$v$. As the interface is traversed from $u$ to $v$, 
the hexagons which border the 
interface on the right are always white and the hexagons along the 
interface on the left are always black.
All the other interfaces will be loops. 

One can imagine generating the interface from $u$ to $v$ 
in the following way. Initially only the boundary hexagons are 
colored. The first step of the walk starts at $u$ and is determined
by the colored boundary hexagons.
When the walk takes a step from $\omega(i-1)$ to $\omega(i)$, 
we refer to the hexagon centered at 
$\omega(i) + [\omega(i)-\omega(i-1)]$ as the hexagon that the 
walk is ``hitting.'' 
If the hexagon the walk is hitting is already colored 
then the color tells us which direction to turn. 
If it is white we turn left, and if it is black we turn right. 
If this hexagon has not been colored yet, we randomly color
it white or black with probability $1/2$. From the point of
view of the walk, we are randomly choosing one of the two 
directions for the next step. Note that hitting a colored hexagon means that 
some other site on the hexagon is occupied. Hence one of the two 
possible next steps will take the walk to a trapping site. 
A little thought reveals that the rule which specifies which direction 
to turn when we encounter a colored hexagon is the same as the rule for 
the SKSAW. The colors simply provide a method for determining which 
direction will keep the endpoint of the walk from being a trapping site.
This relation between the SKSAW and percolation shows that the chordal SKSAW
from $u$ to $v$ has the same distribution as the chordal SKSAW
from $v$ to $u$. Simple examples show that this reversibility 
does not hold for the chordal SKSAW on the square lattice.

A celebrated result of Smirnov is that the scaling limit of 
the percolation interface defined above on the hexagonal lattice
is chordal SLE$_6$ \cite{smirnov,camia_newman}. 
So on the hexagonal lattice the scaling limit of the chordal 
SKSAW is chordal SLE$_6$. Thus it is natural to 
expect (at least for the hexagonal lattice) that the scaling limit
of the full-plane SKSAW is full-plane SLE$_6$, and the scaling limit
of the radial SKSAW is radial SLE$_6$. A sketch of a proof of 
this latter assertion for the hexagonal lattice
has been given by Werner \cite{werner2007lectures}, and a complete proof
has been given by Jiang \cite{jiang}.
If one believes that the scaling limit of the SKSAW does not depend 
on the lattice used, then for all lattices the scaling limit should 
be the appropriate version of SLE$_6$.

On the hexagonal lattice the SKSAW also has a simple interpretation
as a self-avoiding walk with an attractive interaction, and this 
suggests a connection with the tri-critical point or 
$\Theta$-point for polymers. We first review 
the tri-critical point for polymers. It is most easily described
in the context of an interacting self-avoiding walk.
The simplest interacting self-avoiding walk (ISAW) is defined
as follows. For an $n$-step self-avoiding walk $\omega$ let 
\bea
H(\omega) = - \sum_{i=0}^{n-2} \sum_{j=i+2}^n \delta_{|\omega(i)-\omega(j)|,1}
\label{isaw_ham}
\eea
The limits on the sums are such that this energy function counts
the number of pairs of sites on the walk that are nearest neighbors,
not counting those pairs that are only one time step apart. 
The probability of a walk is then taken to be proportional to 
$e^{-\beta H(\omega)}$. The minus sign in \reff{isaw_ham} makes this
an attractive interaction. There is expected to be a phase transition
in $\beta$. For $\beta$ below the critical value the scaling limit
should be the same as for $\beta=0$, the ordinary SAW. 
For $\beta$ above the critical value the scaling limit is expected 
to be a collapsed phase in which the walk is contained in a ball whose
radius is proportional to $N^{1/d}$, where $N$ is the number of sites
in the walk. This critical point is sometimes called the $\Theta$-point. 

We now return to the SKSAW on the hexagonal lattice.
We consider the full-plane SKSAW and run it until there are $N$ steps.
We leave it as an exercise for the reader to check that the 
probability of a walk $\omega$ is proportional to $2^{-N(\omega)}$
where $N(\omega)$ is the number of hexagons that have at least 
one boundary edge that belongs to $\omega$. 
(The easiest way to derive this formula 
is to recognize that each time the walk encounters a hexagon 
that has not been counted in $N(\omega)$ yet, 
the probability acquires a factor of $1/2$. )
So the probability is proportional to $e^{-\beta N(\omega)}$ with 
$\beta=\ln 2$. This weight favors walks that fold back on themselves
often to minimize $N(\omega)$. So it is a sort of 
attractive interaction. 
The energy function $N(\omega)$ is not the same 
as $H(\omega)$ in \reff{isaw_ham}, but it is natural to expect it is 
in the same universality class. It is not obvious that 
$\beta= \ln 2$ is the critical value for $N(\omega)$, but the 
connection with critical percolation suggests that it is. 
Thus we expect that the SKSAW is in the same universality class 
as the tri-critical polymer, and so the scaling limit of the tri-critical
polymer should also be SLE$_6$.

The SKSAW model was introduced in the physics literature under two
different names. Kremer and Lyklema introduced a model 
which they called the infinitely growing self-avoiding walk 
\cite{kremer1985IGSAW}.
It is the full-plane version of the model we have been considering. 
Independently Weinrib and Trugman introduced the model, calling it 
the smart kinetic walk. They considered both the full-plane SKSAW
and the loop SKSAW \cite{weinrib_trugman}.
They showed that the loops 
generated by the loop SKSAW on the hexagonal lattice have exactly the 
same distribution as interfaces in critical percolation on the 
hexagonal lattice. 
The walk is ``smart'' in the sense that it never allows itself to 
become trapped. (Note that the notion of trapping is always with 
respect to the site the walk is trying to reach.) The walk is 
``kinetic'' in the sense that it can be defined by a simple 
rule that gives the transition probabilities at each step.
We refer to the model as the smart kinetic self-avoiding walk
since in all the versions the walk is 
smart and kinetic, while it is infinitely growing only 
in some of the versions. 

Several other random walk models with some form of self-avoidance 
were introduced at about the same time.
Amit, Parisi, and Peliti defined what they called the true self 
avoiding walk (TSAW) \cite{amit_et_al}. 
It is a nearest neighbor random walk in 
which the walk is allowed to visit sites it has visited before, 
but is discouraged from doing so. More precisely, the walk can jump only 
to nearest neighbors and the probability of jumping to the nearest neighbor
$u$ is proportional to $e^{-g n_u}$ where $n_u$ is the number of previous
visits to $u$, and $g>0$ is a parameter. 
They argued that this model is in a different universality class
from the self-avoiding walk and that the upper critical dimension is two.
(In their paper they often refer to the TSAW as simply the self-avoiding
walk and use the term self-repelling polymer-chain for what is now 
called the self-avoiding walk.) 
Another random walk model is the Kinetic Growth Walk (KGW)
\cite{majid1984kinetic}. It randomly chooses one of its 
unoccupied nearest neighbors at each step. Unlike the SKSAW it is 
allowed to move to trapping sites, so the walk can get trapped.
On certain oriented lattices the KGW cannot get trapped and so is equivalent 
to the SKSAW \cite{manna1989kinetic}.
Neither the TSAW or the KGW are 
expected to be in the same universality class as the SKSAW.  
The SKSAW is the same as the Laplacian random walk with the parameter 
set to zero. The Laplacian random walk was introduced in 
\cite{lyklema1986laplacian}, and its connection with SLE was studied in 
\cite{lawler2006laplacian}.

As noted earlier, the connection of the SKSAW on the hexagonal 
lattice with critical percolation was studied in \cite{weinrib_trugman}.
Ziff, Cummings, and Stell introduced a random walk that generates 
the perimeter of percolation clusters for critical site percolation
on the square lattice \cite{ziff1984}. 
A random walk model that generates the perimeters for bond percolation 
on the square lattice was given by Gunn and Ortu\~no \cite{gunn_ortuno}.

Coniglio, Jan, Majid, and Stanley \cite{coniglio1987} 
studied the SKSAW on the hexagonal lattice and argued that it 
is a polymer model at the tri-critical point. 
Since the SKSAW on the hexagonal lattice is equivalent to a critical 
percolation interface, they concluded that it appears 
that this particular polymer chain has the same statistics as
percolation interfaces.
Duplantier and Saleur \cite{duplantier_saleur}
studied a model of a self-avoiding walk 
on the hexagon lattice with random forbidden hexagons. 
They argued that this model describes the 
tri-critical point of the self-avoiding walk. They then used Coulomb gas 
methods to derive exact values for the critical exponents. 

Since the scaling limit of percolation is SLE$_6$, the scaling limit
of the tri-critical point for polymers should be described by SLE$_6$.  
Gherardi studied this prediction for the ISAW with Monte Carlo 
simulations of the ISAW and SLE$_6$ \cite{gherardi}.
There are some subtleties related to just 
how we should define the ISAW in a bounded domain to obtain 
a scaling limit of radial or chordal SLE$_6$ that we discuss in 
the final section of this paper.

In this paper we point out some consequences of the conjecture
that the full-plane SKSAW converges to full-plane SLE$_6$
and use these observations to test the conjecture with 
simulations of the SKSAW on the square lattice.
We perform the same simulations on the hexagonal lattice for 
comparison.
In the next section we review full-plane SLE$_6$ and state the 
conjecture for the scaling limit of the SKSAW 
and its consequences. In sections three and four 
we present the results of our simulations.
Our conclusions and some open questions are 
discussed in the final section.

\section{The conjecture and its consequences}

We begin this section by reviewing the definition of full-plane SLE
following section 6.6 of \cite{lawler_book}. 
The Loewner equation for full-plane $SLE$ is 
\bea
{\dot g}_t(z) = g_t(z) \frac{e^{-iU_t}+g_t(z)}{e^{-iU_t}-g_t(z)}
\eea
with the initial condition $\lim_{t \ra -\infty} e^t g_t(z) = z$
for $ z \in \complex \setminus \{0\}$. 
The driving function $U_t$ is defined as follows.
Let $B^1_t$ and $B^2_t$ be independent Brownian motions with mean zero
and $E (B^i_t)^2 =t$. Let $Y$ be another independent random variable
which is uniformly distributed on $[0,2\pi]$. 
For $t \ge 0$, $U_t = \sqrt{\kappa} B^1_t + Y$ and for 
$t \le 0$, $U_t = \sqrt{\kappa} B^2_{-t} + Y$.
For each $t$, $g_t(z)$ is a conformal map of a domain $H_t$ onto 
$\complex \setminus \overline{\disc}$. There is a curve $\gamma(t)$
defined for $t \in (-\infty,\infty)$ such that $H_t$
is the hull of $\gamma[-\infty,t]$ (the unbounded component of 
$\complex \setminus \gamma[-\infty,t]$).
$\gamma(-\infty)$ is defined to be $0$. 
One can think of full-plane SLE as the limit as 
$\epsilon \ra 0$ of radial SLE in the region $|z| > \epsilon$
where the radial SLE starts at a random point on the boundary of 
the small disc and ends at $\infty$. 

We now let $\gamma$ denote a full-plane SLE$_6$ starting at the origin, 
and we let $W(t)$ be a two dimensional Brownian motion starting 
at the origin. 
Consider a simply connected domain $D$ (other than the full plane)
containing the origin.  
Let $\sigma_D$ be the time when $\gamma(t)$ first exits this domain, and 
$\tau_D$ the time when $W(t)$ first exits the domain. 
So $\gamma(\sigma_D)$ and $W(\tau_D)$ are random points on the boundary 
of $D$. The distribution of the latter is of course harmonic
measure. The curve $\gamma$ will intersect itself without crossing 
itself. The hull of $\gamma[0,t]$ is defined to be 
the curve up to time $t$ 
along with the regions that are enclosed by the curve. More precisely,
it is the complement of the unbounded component of the complement of 
$\gamma[0,t]$. We denote the hull of $\gamma[0,t]$ by $K_t$. 
The hull of $W([0,t])$ is defined similarly and denoted by $\hat{K}_t$. 
There is a deep connection between these two hulls.
 
\begin{proposition}(Propositions 6.31 and 6.32 in \cite{lawler_book})

\noindent 
(i) $\gamma(\sigma_D)$ and $W(\tau_D)$ have the same distribution, i.e., 
harmonic measure. \\

\noindent
(ii) The hulls $K_{\sigma_D}$ and $\hat{K}_{\tau_D}$ have the same
distribution.
\end{proposition}

Let $v \in \partial D$. If we condition $\gamma(t)$ on the event 
$\{\gamma(\sigma_D)=v\}$ and condition $W(t)$ on the event
$\{W(\tau_D)=v\}$ then we have two probability measures on hulls
in $D$ that intersect the boundary at $v$. The proposition implies that
the two probability measures are the same. 
This probability measure is an example of a radial restriction measure. 
Radial restriction measures were studied by Wu \cite{wu}.
We refer the reader there for the full definition. Informally, 
the restriction property is defined as follows. 
Given a simply connected domain $D$ and points $u$ on the boundary and 
$v$ in the interior we say that a closed subset $K$ of $\overline{D}$
is a hull if $K$ is connected, $\complex \setminus K$
is connected, $v \in K$ and $K \cap \partial D = \{u\}$. 
Suppose that for each simply connected domain $D$ (other than $\complex$)
and points $u$ on the boundary and $v$ in the interior 
we have a probability measure $P_{D,u,v}$ on hulls. We assume that 
these measures are related to each other by conformal transformations
in the obvious way.
(So they are completely determined by the measure for a single 
choice of $D$, $u$, $v$.) We say that these measures satisfy the 
restriction property if whenever $D^\prime \subset D$ with $u$ on the 
boundaries of $D^\prime$ and $D$ and $v \in D^\prime$, then $P_{D^\prime,u,v}$ is 
just $P_{D,u,v}$ conditioned on the event that $K \subset \overline{D^\prime}$. 

The conformal invariance of Brownian motion implies that 
the measure in part (ii) of the proposition is conformally 
invariant. This invariance also allows us to do explicit 
calculations for the measure. Take $D$ to be the unit disc $\disc$ and 
condition on the event that we exit $D$ at the point $1$. 
Then it is an easy exercise to show that 
\bea
P(K \cap A = \emptyset) = \Phi_A^\prime(1)
\label{radial_restrict_bm}
\eea
where $A$ is any closed subset of $\overline{\disc}$ such that 
$A =\overline{ \disc \cap A}$, $\disc \setminus A$ is simply connected,
the origin is in $\disc \setminus A$, and $1$ belongs to its boundary.
$\Phi_A$ is the conformal map of $\disc \setminus A$ onto $\disc$ that
fixes the origin and $1$. 
The formula holds for all hulls $K$ in $\overline{\disc}$ with $0 \in K$
and $K \cap \partial \disc = \{1\}$. 
In general, radial restriction measures are completely 
characterized by a generalization of this formula:
\bea
P(K \cap A = \emptyset) = |\Phi^\prime_A(0)|^\alpha \, \Phi_A^\prime(1)^\beta
\label{radial_restrict}
\eea
where the real parameters $(\alpha,\beta)$ completely determine
the restriction measure. So our measure is the particular case of 
$\alpha=0, \beta=1$.

We now state the main conjecture studied in this paper and 
the consequences that would follow from the conjecture 
using the rigorous results described above. 
For the hexagonal lattice these are all theorems.
This conjecture is implicit in the literature; in particular 
the analogous conjecture for the half-plane may be found in 
chapter 1 of \cite{lawler_book}.

\medskip
\noindent {\bf Conjecture:}
{\it In the limit that the lattice spacing converges to zero, the 
full-plane smart kinetic self-avoiding walk from $0$ to $\infty$ on any regular 
lattice converges in 
distribution to full-plane SLE$_6$ from $0$ to $\infty$.}

\medskip
\noindent {\bf Consequence 1:}
{\it Let $D$ be a simply connected domain containing the origin.
We refer to the distribution of the point where the 
smart kinetic self-avoiding walk first exits the domain $D$ as the 
exit distribution for $D$.
As the lattice spacing converges to zero the exit distribution
converges to harmonic measure. }

\medskip
\noindent {\bf Consequence 2:}
{\it Let $D$ be a simply connected domain containing the origin, and 
let $v$ be a point on its boundary.
We consider the hull of the smart kinetic self-avoiding walk up to the time it 
exits $D$, conditioned on the event that it exits at $v$.
In the limit that the lattice spacing goes to zero this process 
is conformally invariant. For $D$ equal to the unit disc and $v=1$, 
the scaling limit satisfies \reff{radial_restrict_bm}.}

\section{Simulations}

We generate samples of the SKSAW by growing the walk according 
to the definition given in the first section. 
Some details of how we do this are given at the end of this
section.
We test the conjecture that the scaling limit of the full-plane 
SKSAW is full-plane SLE$_6$ for the square lattice.
In this section we discuss two types of tests that assume the 
scaling limit is rotationally invariant. 
We perform the same simulations for the hexagonal lattice 
for comparison purposes.
In the next section we discuss two similar types of tests that do not 
assume the scaling limit is rotationally invariant. 

For the first type of simulation we compute the exit distribution
for three particular domains and compare the results with harmonic
measure for those domains. 
The first domain is a disc with radius $2$ centered at the point $1$:
\beann
D_1=\{ z : |z-1|<2\}.
\eeann
The second domain is a horizontal strip of width $3$ with the distance
from the origin to the top boundary equal to $2$. So 
\beann
D_2 = \{ z : -1 < Im(z) < 2 \}
\eeann
The third domain is an equilateral triangle with the origin at the 
center.
\beann
D_3 = hull\{ (2,0), (-1,\sqrt{3}), (-1,-\sqrt{3})\}
\eeann
where $hull$ means the convex hull of the three points.
All three domains are scaled so that the distance from the origin 
to the boundary of the domain is $1$. 

We are interested in the limit that the lattice spacing goes to zero.
One effect of the non-zero lattice spacing is that the point at 
which the walk exits the domain is a discrete random variable. 
If the scaling limit is rotationally invariant then it does not 
matter how the lattice is oriented with respect to the domain.
So if we assume the scaling limit is rotationally 
invariant, then we can reduce the effect of this discreteness 
by averaging over rotations of the lattice.
We make this assumption and implement this averaging 
by randomly choosing a rotation of 
the domain for each sample that we generate. 

The second type of simulation test in this section
is based on the observation that 
if the conjecture is true, then if we condition on where the process
exits the domain $D$ we have a radial restriction measure whose hull has
the same distribution as Brownian motion conditioned to exit the 
domain at the same point. 
We consider the process in the unit disc started at the origin. 
We would like to condition on the event that 
it exits at $1$. Assuming that the scaling limit is rotationally 
invariant, in the scaling limit this is equivalent to applying a random 
rotation to take the endpoint of the process to $1$. 
This gives a curve in the unit disc from $0$ to $1$.
Eq. \reff{radial_restrict_bm} can be used 
to explicitly calculate the distribution of various random variables
for this process. For our simulations we consider three random variables.
The random variable $X$ is the farthest the curve travels to 
the left. The random variable $Y$ is the farthest the curve travels
in the upwards direction, and $Z$ is the maximum distance the curve
travels from $1$. So letting $\omega(t)$ denote the curve, 
\bea
X &=& \max_t  \,  \left[ - Re(\omega(t)) \right] \nonumber \\
Y &=& \max_t \, Im(\omega(t)) \nonumber \\
Z &=& \max_t \, |\omega(t) -1| \nonumber \\
\label{xyzdef}
\eea
Eq. \reff{radial_restrict_bm} can be used 
to compute the distribution of these three random variables.
We omit the details of these computations. 

For all our simulations we display the results using the cumulative
distribution function (CDF) rather than the density. 
(Computing the density from a simulation requires taking a numerical 
derivative and so adds further uncertainty.)
For the exit distributions for the three domains, we let $\Theta$ be the 
polar angle of the point where the walk first exits the domain and 
look at the CDF defined by $F_\delta(\theta)=P_\delta(\Theta \le \theta)$.
The subscript $\delta$ is the lattice spacing. We let 
$H(\theta)$ be this CDF for harmonic measure. We are testing the 
conjecture that $F_\delta(\theta)$ converges to $H(\theta)$ 
as $\delta \ra 0$. If we plot these two functions together they 
are indistinguishable. So we will only show plots of the difference 
$F_\delta(\theta)- H(\theta)$. 

\begin{figure}[tbh]
\includegraphics{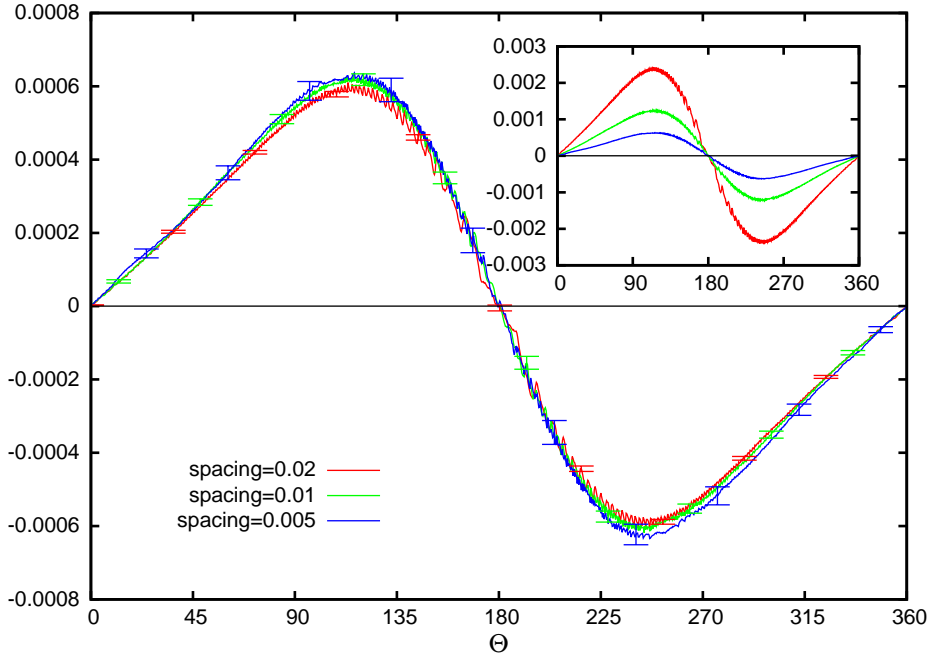}
\caption{\leftskip=25 pt \rightskip= 25 pt 
The differences between the simulation of the CDF for the exit 
distribution for the SKSAW on the square lattice
and the CDF for harmonic measure for domain $D_1$ are shown in the inset. 
In the main plot these differences are rescaled: 
$\delta=0.005$ is unchanged while $\delta=0.01$ and 
$0.02$ are rescaled by factors of $1/2$ and $1/4$, respectively. 
}
\label{fig_sq_off_center}
\end{figure}

We have run simulations for both the square lattice and the hexagonal
lattice. For the hexagonal lattice the conjecture that the
scaling limit is full-plane SLE$_6$ is a theorem.
The purpose of these simulations is simply to provide 
some comparison with the square lattice simulations.
So we focus our description of the results on the square lattice. 

For all three domains we have run simulations with 
$\delta=0.02, 0.01, 0.005$. 
For each lattice spacing we generated $10^9$ samples.
For the square lattice the differences for the 
domain $D_1$ (a disc whose center is not at the origin) 
are shown in the inset in figure \ref{fig_sq_off_center}.
This plot clearly shows that the differences are going to zero
as $\delta \ra 0$. It
also suggest a sort of finite size scaling with the magnitude of 
the differences proportional to $\delta$.
In the larger plot in figure \ref{fig_sq_off_center}
we test this scaling as following. 
The difference for $\delta=0.005$ is plotted unchanged. The difference
for $\delta=0.01$ is multiplied by $1/2$, and the difference for 
$\delta=0.02$ is multiplied by $1/4$. 
If the differences are indeed proportional to $\delta$,
then this rescaling should cause the curves to collapse to 
a single curve. In the plot of the rescaled differences we see that
they collapse together so well that 
it is difficult to distinguish the three curves.
The error bars shown represent plus or minus two standard deviations
for the statistical errors. 

\begin{figure}[tbh]
\includegraphics{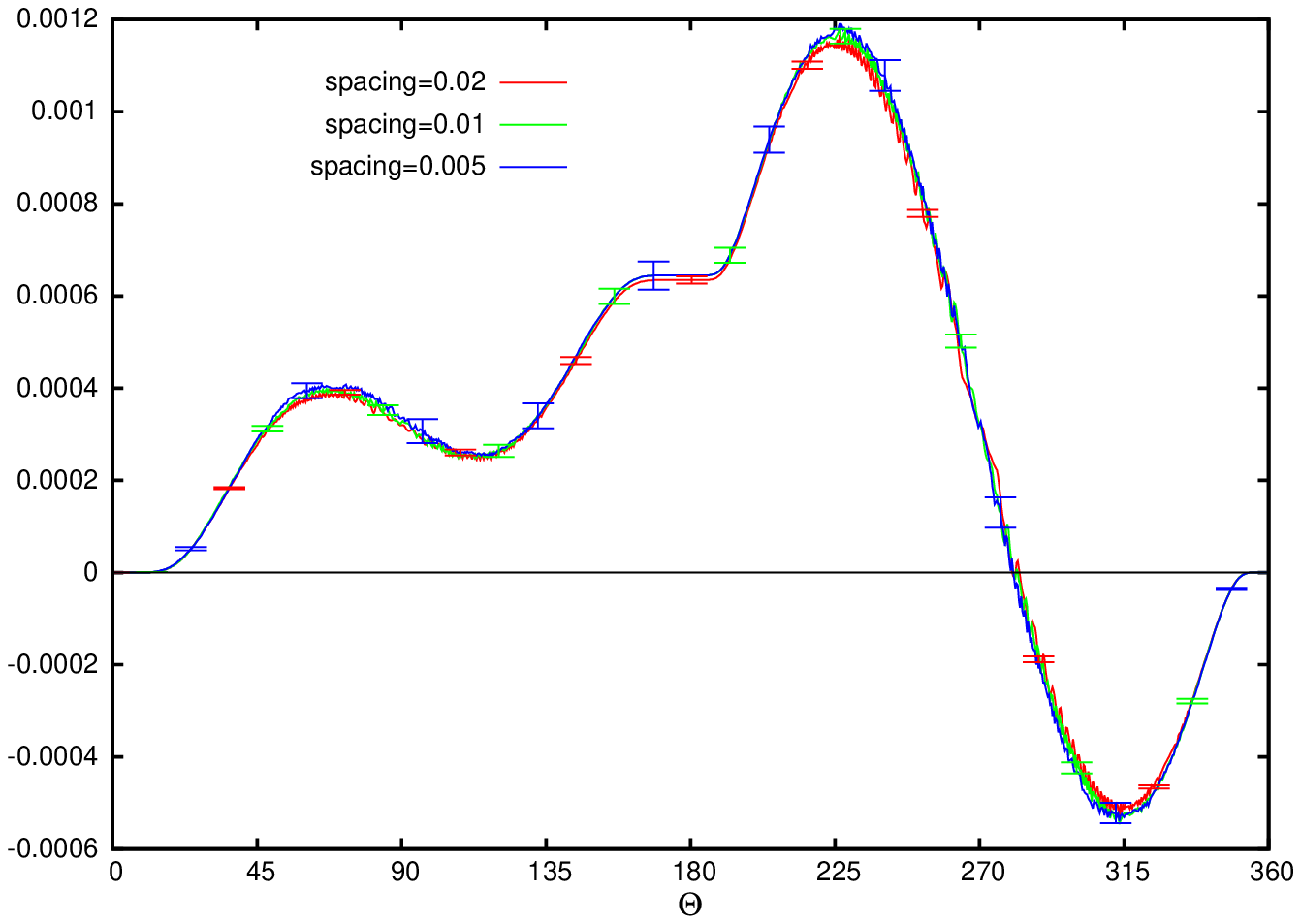}
\caption{\leftskip=25 pt \rightskip= 25 pt 
Differences for the square lattice for the domain $D_2$. 
They are rescaled as we did for $D_1$.
}
\label{fig_sq_strip}
\end{figure}

\begin{figure}[tbh]
\includegraphics{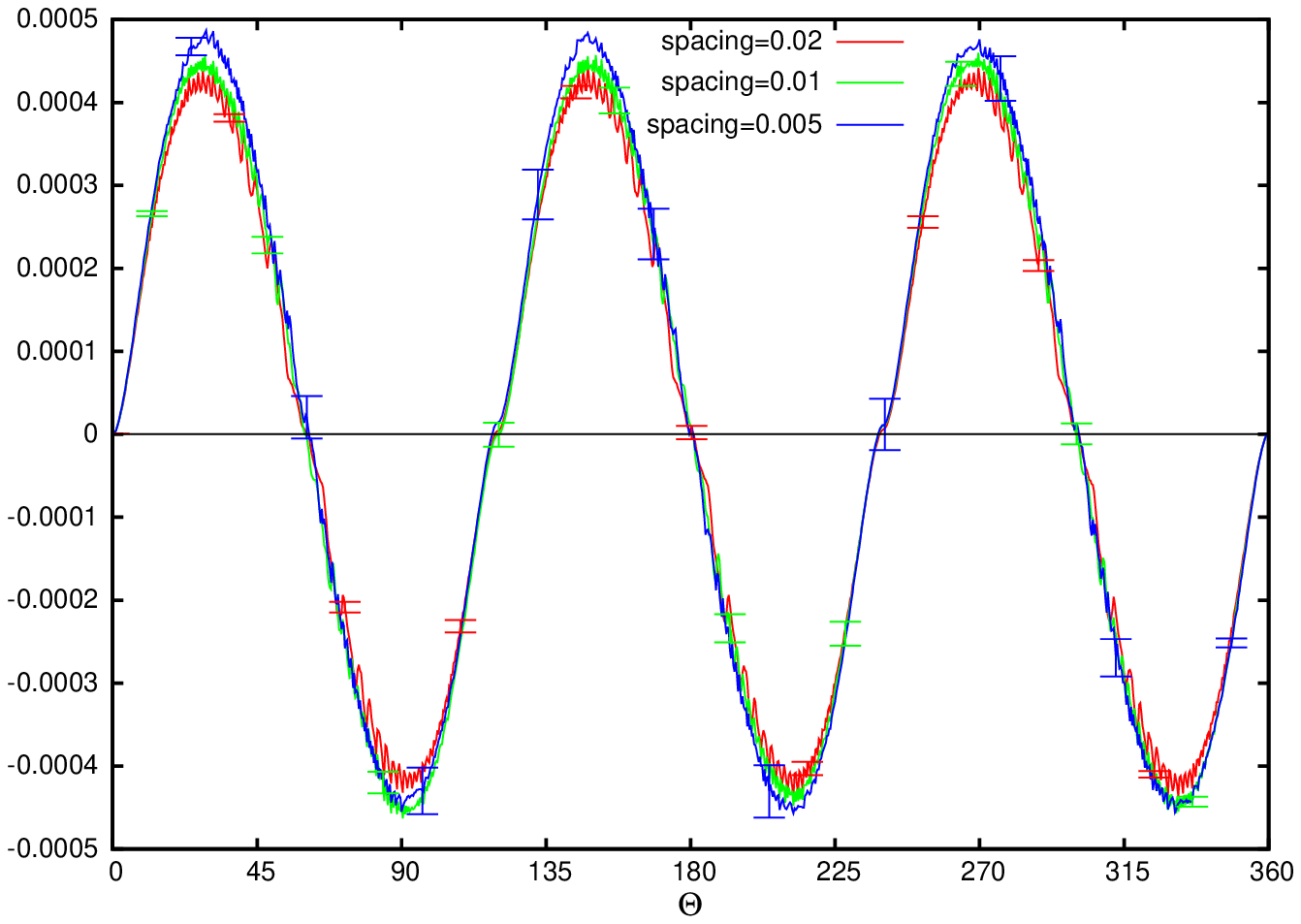}
\caption{\leftskip=25 pt \rightskip= 25 pt 
Rescaled differences for the square lattice for domain $D_3$.
}
\label{fig_sq_triangle}
\end{figure}

In figure \ref{fig_sq_strip} we plot the differences for the square lattice
for the domain $D_2$ (the strip). The differences are 
rescaled as we did for the previous domain. ($\delta=0.005$ is unchanged,
$\delta=0.01, 0.02$ are multiplied by $1/2, 1/4$ respectively.)
Again we see that the rescaled differences collapse nicely onto a 
single curve. 
Finally, for the square lattice and domain $D_3$ (the triangle), 
figure \ref{fig_sq_triangle} shows the rescaled differences. 

\begin{figure}[tbh]
\includegraphics{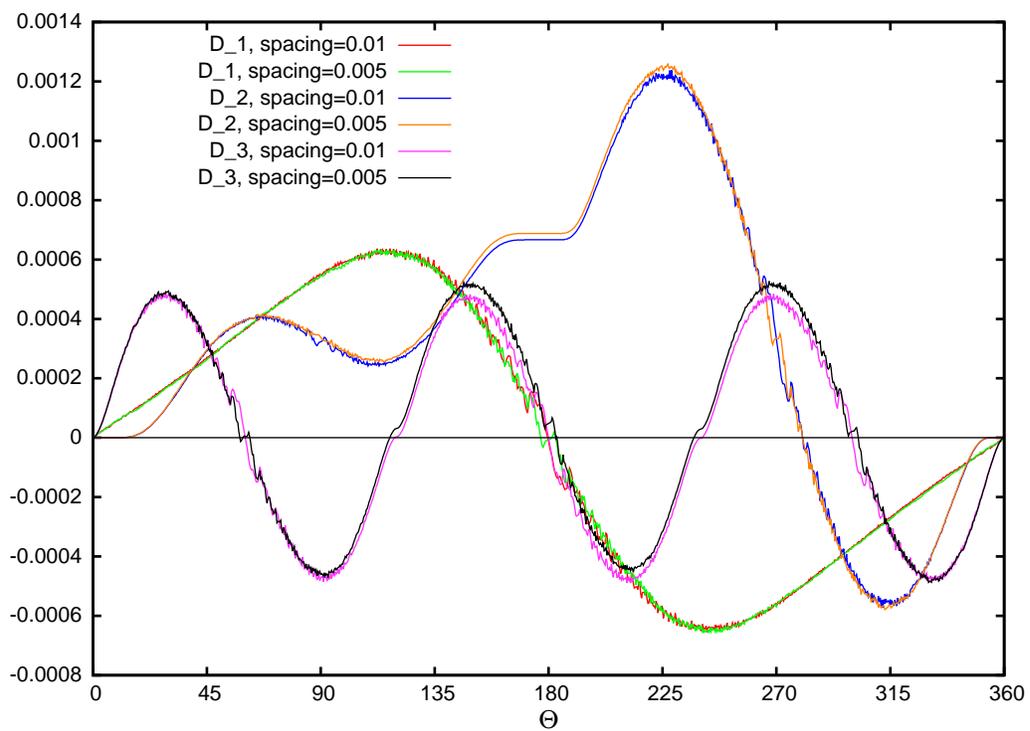}
\caption{\leftskip=25 pt \rightskip= 25 pt 
Differences for the hitting distribution for the hexagonal lattice
for all three domains.
}
\label{fig_hex_all_domain}
\end{figure}

The results for the hexagonal lattice are very similar. The size of 
the difference between the SKSAW CDF and the harmonic measure CDF appears to 
be proportional to $\delta$. When we rescale these differences 
by factors proportional to $1/\delta$, the curves again collapse 
nicely onto a single curve. We show these rescaled curves for the hexagonal
lattice for $\delta=0.01,0.02$ for all three domains in 
figure \ref{fig_hex_all_domain}.
We were quite surprised to find that the shape of these 
curves is similar to the corresponding curves for the 
square lattice. (Their magnitudes are slightly different.)
This is surprising since one would expect 
the difference between the exit distribution on a lattice and 
harmonic measure to depend very strongly on the lattice. 
We should remind the reader that we are assuming the scaling limit is 
rotationally invariant and using this assumption to average 
our simulations over rotations of the lattice with respect to the domain. 
Nonetheless, even with this averaging over rotations we would still 
expect the differences being plotted to depend on the lattice used.

\clearpage

\begin{figure}[tbh]
\includegraphics{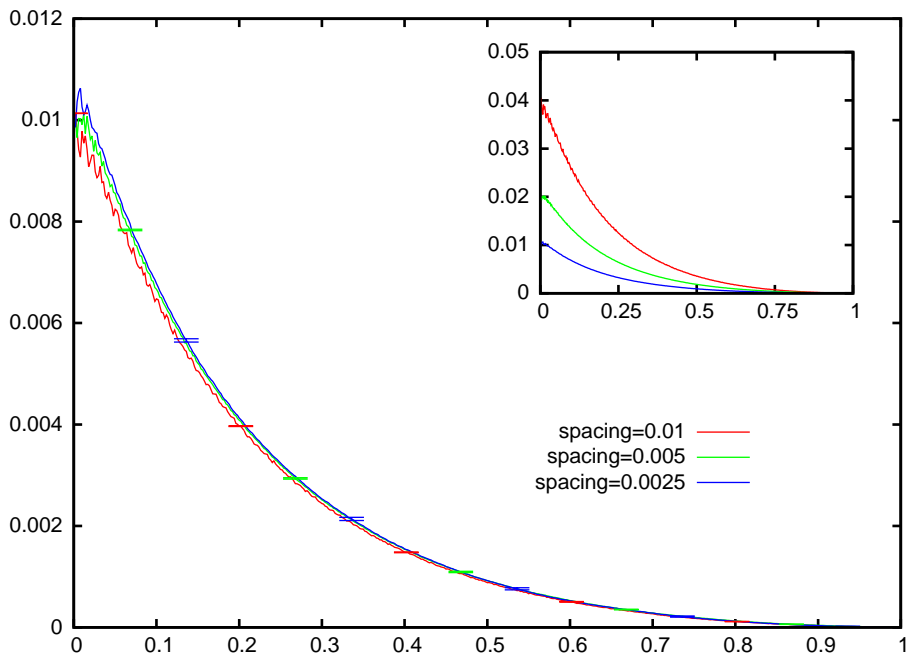}
\caption{\leftskip=25 pt \rightskip= 25 pt 
Differences for the RV $X$ in the disc on the square lattice are 
shown in the inset. The main plot shows these differences rescaled 
by a factor proportional to $1/\delta$.
}
\label{fig_sq_disc_rv0}
\end{figure}

\begin{figure}[tbh]
\includegraphics{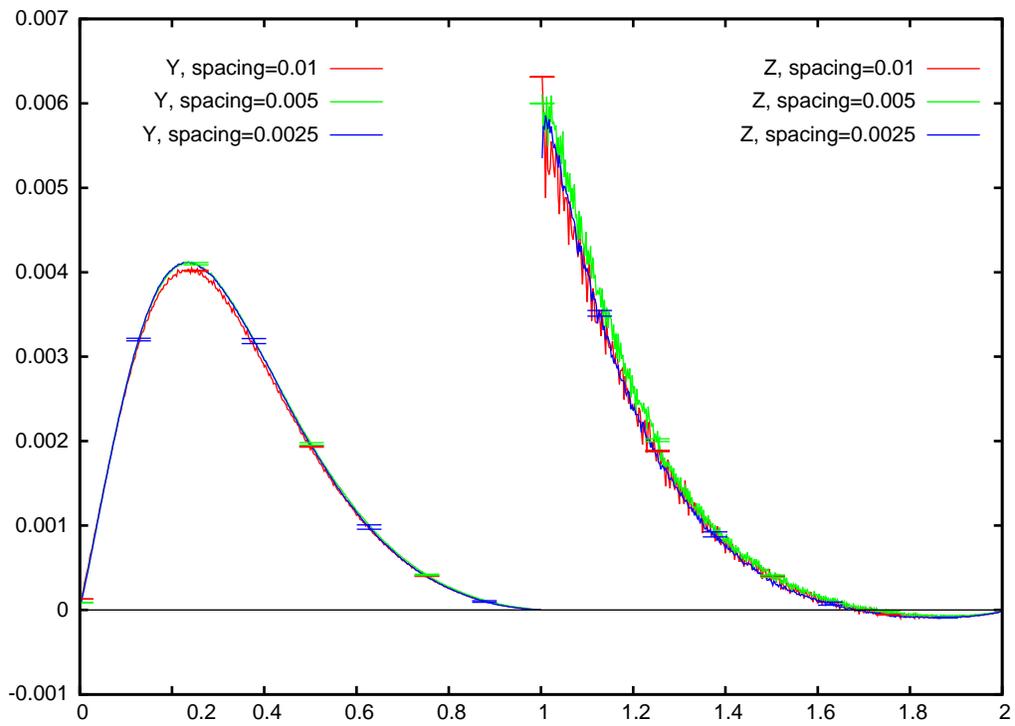}
\caption{\leftskip=25 pt \rightskip= 25 pt 
Rescaled differences for the RV's $Y$ and $Z$ 
in the disc on the square lattice.
}
\label{fig_sq_disc_rv12}
\end{figure}

\begin{figure}[tbh]
\includegraphics{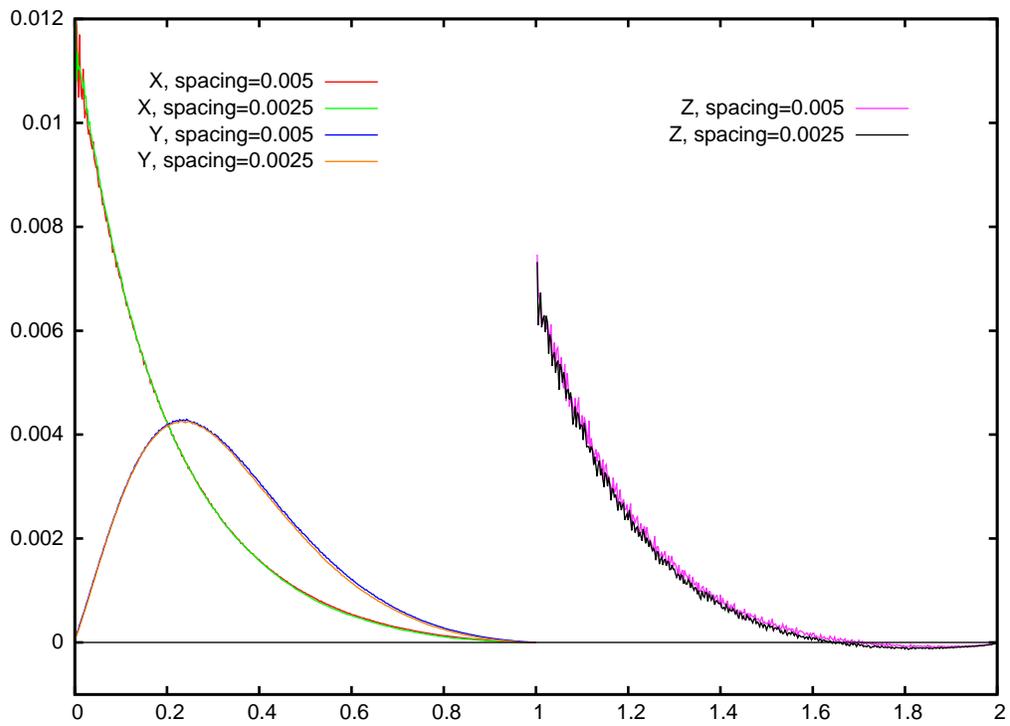}
\caption{\leftskip=25 pt \rightskip= 25 pt 
Differences for all three RV's for the hexagonal lattice.
}
\label{fig_hex_all_rv}
\end{figure}

Our second set of tests of the conjecture that the scaling limit 
of the full-plane SKSAW is full-plane SLE$_6$ is to simulate the 
distributions of the three random variables $X,Y$ and $Z$
defined in \reff{xyzdef}. 
We have performed simulations for both the square and hexagonal lattices
with lattice spacings of $\delta=0.01, 0.005, 0.0025$. 
For each lattice spacing we generated $10^9$ samples.
For each random variable and lattice we plot 
the difference of the CDF from the simulation and the 
exact CDF predicted by SLE$_6$. The inset in 
figure \ref{fig_sq_disc_rv0} shows these 
differences for $X$ for the square lattice. As with the exit 
distributions it appears the size of the difference is proportional
to $\delta$. The main plot in 
figure \ref{fig_sq_disc_rv0} tests this by 
plotting the differences for $\delta=0.005$ 
and $\delta=0.0025$ multiplied by factors of $2$ and $4$, respectively, 
while the difference for $\delta=0.01$ is unchanged.
The three curves collapse 
to a single curve, supporting the hypothesis that the difference
is proportional to $\delta$. Similar rescaled plots of the 
difference for $Y$ and $Z$ for the square lattice are shown in 
figure \ref{fig_sq_disc_rv12}.
The rescaled differences for all three random variables for the hexagonal 
lattice are shown in figure \ref{fig_hex_all_rv}.
Only the differences for $\delta=0.05$ and  $\delta=0.0025$ are shown
in this figure.

\clearpage

In choosing the lattice spacing to use in these simulations
there is a trade-off. One would like to make the lattice spacing 
small so that the simulation is close to the scaling limit. 
But a smaller lattice spacing means the time to generate a sample is 
longer and so fewer samples can be generated resulting in larger 
statistical errors.  We have chosen
lattice spacings that result in relatively short walks.
For the SKSAW in the unit disc with the smallest lattice spacing of
$\delta=0.0025$, the average number of steps is approximately
$19,000$ for the square lattice and $15,800$ for the hexagonal lattice.
For the lattice spacings that we use, we are able to generate $10^9$ samples.
The result is that we can study the deviation of the simulation 
from the scaling limit accurately and see clearly that 
the leading order correction is proportional to $\delta$.
We expect the critical exponent $\nu$ to be $4/7$ \cite{duplantier_saleur},
so the average 
number of steps in the SKSAW should be proportional to $\delta^{7/4}$. 
Thus to reduce the lattice spacing by a factor of $2$ increases
the time required by more than a factor of $3$.
Our simulations for the smallest lattice spacings 
typically took about a week running on 100 CPU's. 

We conclude this section with some details about how we
generate samples of the SKSAW. At each step we need to 
determine which nearest neighbors are occupied sites 
and which are trapping sites. 
We use a hash table to record which sites have been visited 
so that we can efficiently check if nearest neighbors are 
occupied. Determining if they are trapping sites is a bit more 
subtle. 

\begin{figure}[tbh]
\includegraphics{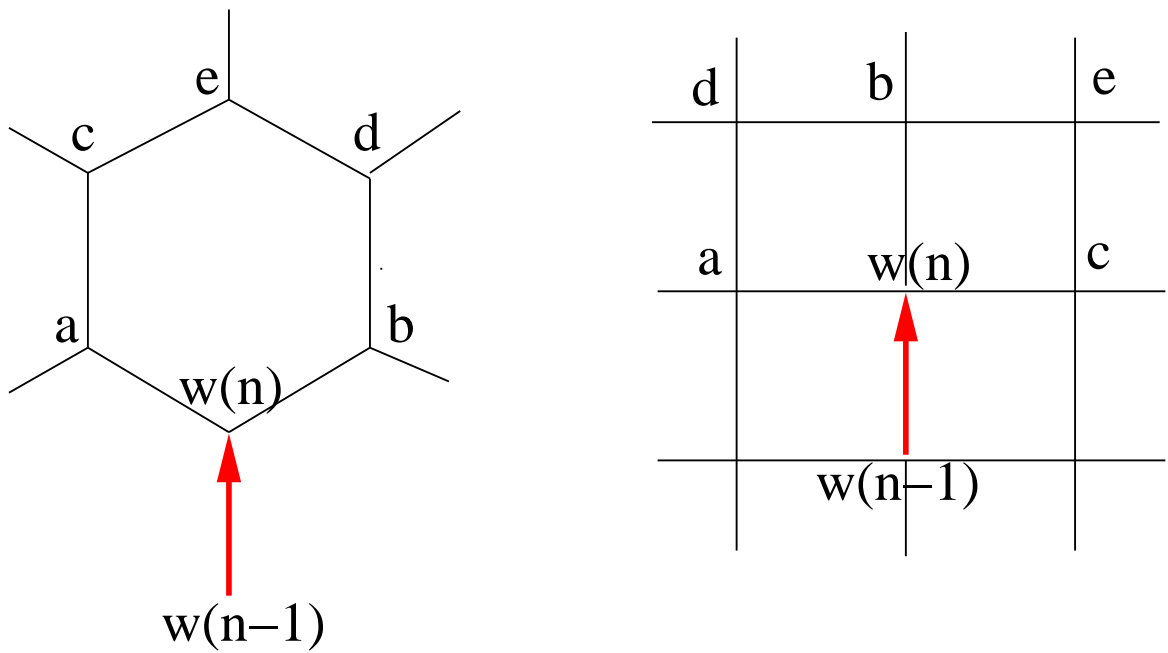}
\caption{\leftskip=25 pt \rightskip= 25 pt 
On the hexagonal lattice, one of the nearest neighbors $a$ and $b$ 
can be a trapping site only if at least one of the sites $c,d,e$
is occupied. 
On the square lattice, one of the nearest neighbors $a,b$ and $c$ 
can be a trapping site only if at least one of the sites $d,b,e$
is occupied. 
}
\label{algorithm}
\end{figure}

Suppose the walk $\omega(i)$ has been defined for $i \le n$,
and let $p$ a nearest neighbor of $\omega(n)$. 
Since the current site $\omega(n)$ is not a trapping site, 
in order for $p$ to be a trapping site the walk must be close to 
forming a loop that disconnects $p$ from $\infty$. 
To determine which sites are trapping sites we need to know 
which sites are inside and which are outside the loop that is 
about to be formed. This can be efficiently determined using the 
winding angle of the loop \cite{kremer1985IGSAW}.
We first consider the hexagonal lattice, referring to figure \ref{algorithm}.
If one of the nearest neighbors $a$ or $b$ is occupied then 
the other nearest neighbor cannot be occupied and cannot be a trapping site. 
So we need only consider the case that neither nearest neighbor
is occupied. We describe how we determine if the nearest neighbor 
$a$ is a trapping site. The procedure for $b$ is similar. 
Since $\omega(n)$ is not a trapping site, 
in order for $a$ to be a trapping site, 
at least one of the sites $c, d, e$ must be occupied. The portion of the walk 
from this occupied site to $\omega(n)$ is then close to forming a loop.
We can close it to form a loop by adding two or three steps. 
We do this in such a way that the extension does not pass 
through $a$. We then compute the orientation
of the loop we have formed.
This orientation tells us if $a$ is inside or outside the 
loop, and so we can determine if the loop disconnects the nearest 
neighbor $a$ from $\infty$. If it does then the nearest neighbor $a$
is a trapping site. (We note that the process we have just described to 
form a loop is only to determine which sites are trapping sites. 
The two or three steps added to form this loop should not be 
confused with the next step that the walk will take.) 
Once we have determined which nearest nearest neighbors of 
$\omega(n)$ are not allowable since they occupied or trapping, we then
pick one of the allowable nearest neighbors with equal probability
to be $\omega(n+1)$. 

For the square lattice we refer again to figure \ref{algorithm}.
Since $\omega(n)$ is not a trapping site, 
in order for one of the nearest neighbors $a, b, c$ 
of $\omega(n)$ to be a trapping site, at least one of the sites 
$a, b, c, d, e$ must be occupied. 
For example, if $b$ is occupied then one of $a$ or $c$ must be a trapping
site. (Which one depends on how the walk connects $b$ to $\omega(n)$.) 
If $a$ is occupied then it is possible to connect $a$ to $\omega(n)$ 
in such a way that $b$ and $c$ are trapping sites. However, this 
connection would mean that $\omega(n)$ is a trapping site. So 
this case does not happen. Likewise we need not worry about the 
possibility that $c$ being occupied makes $a$ and $b$ trapping sites.
If $d$ is occupied then depending on how the walk connects $d$ to 
$\omega(n)$ either $a$ will be a trapping site or $b$ and $c$ will
both be trapping sites. Similarly, if $e$ is occupied then 
either $c$ is a trapping site or both $a$ and $b$ are trapping sites. 
When one of $b,d,e$ are occupied, we determine which nearest neighbors are 
forced to be trapping sites by extending the walk one or two steps 
from $\omega(n)$ to the occupied site in such a way that the walk 
does not pass through the nearest neighbor in question. 
The orientation of this loop then determines whether the site is inside the 
loop. If it is inside it is a trapping site.

To compute the orientation of loops efficiently, as we generate the 
walk we record the total of the angles of the turns the walk has 
made up to the current point. 
Given a loop we can compute the total of the angles of the turns in the loop
by taking the difference of two of these totals. This total for the 
loop tells us the orientation of the loop.
In principle the time needed to generate
a sample with $N$ steps should be $O(N)$, but 
we have not tested if this is true in practice.

\section{Further simulations}

In the previous section we assumed that the scaling limit of the 
SKSAW is rotationally invariant. When we studied the exit 
distribution for various domains, we used this assumption to 
average over rotations of the lattice with respect to the domain. 
This averaging smooths out the discrete nature of the exit distribution
and allows for a better comparison with harmonic measure. 
For the random variables $X,Y,Z$ in the disc, we used this rotational
invariance to rotate the lattice so that the SKSAW exits the disc at 
$1$. In this section we do not assume that the scaling limit is 
rotationally invariant. We carry out two types of simulations that are
similar to those of the previous section. All the simulations in 
this section were only done for the square lattice. 

\begin{figure}[tbh]
\includegraphics{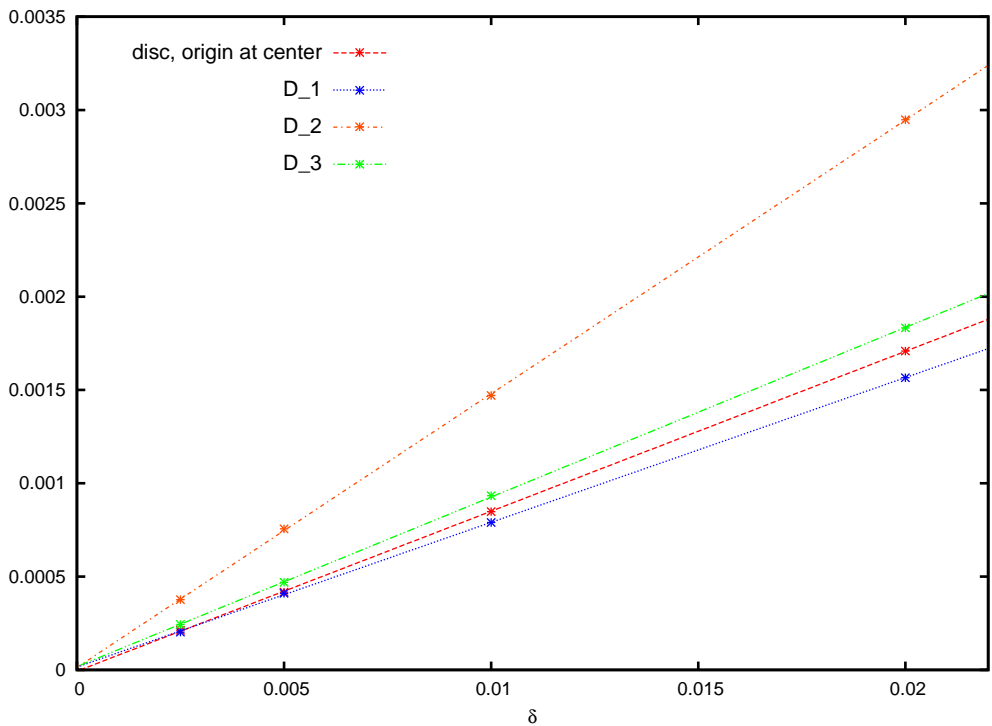}
\caption{\leftskip=25 pt \rightskip= 25 pt 
The $L^1$ norms of the differences between the simulation of the CDF 
for the exit distribution for the SKSAW on the square lattice 
and the CDF for harmonic measure for the four domains.
In these simulations the orientation of the lattice with respect to the 
domain is kept fixed. 
}
\label{fig_rt_exit}
\end{figure}

The first set of simulations compute the exit distribution for the 
same three domains considered in the previous section, 
but without averaging over the orientation of the lattice. 
We also compute the exit distribution for a unit disc when the SKSAW starts
at the center of the disc and the orientation of the lattice is fixed. 
If we were to average over orientations of the lattice,
then this distribution would trivially be uniform on the boundary. 
Without the averaging it is not. Checking that it converges to the uniform
distribution as the lattice spacing goes to zero 
serves as a test of the rotation invariance of the scaling limit. 
For each of the four domains we simulated the SKSAW with lattice 
spacings of $\delta=0.02, 0.01, 0.005, 0.0025$. 
For each lattice spacing we generated $10^9$ samples.
We computed the $L^1$ norm of the difference between the CDF of the 
exit distributions from the simulations and the CDF of harmonic measure. 
These norms are plotted as a function of $\delta$ in figure 
\ref{fig_rt_exit}. The lines shown are least squares fits to the data. 

\begin{figure}[tbh]
\includegraphics{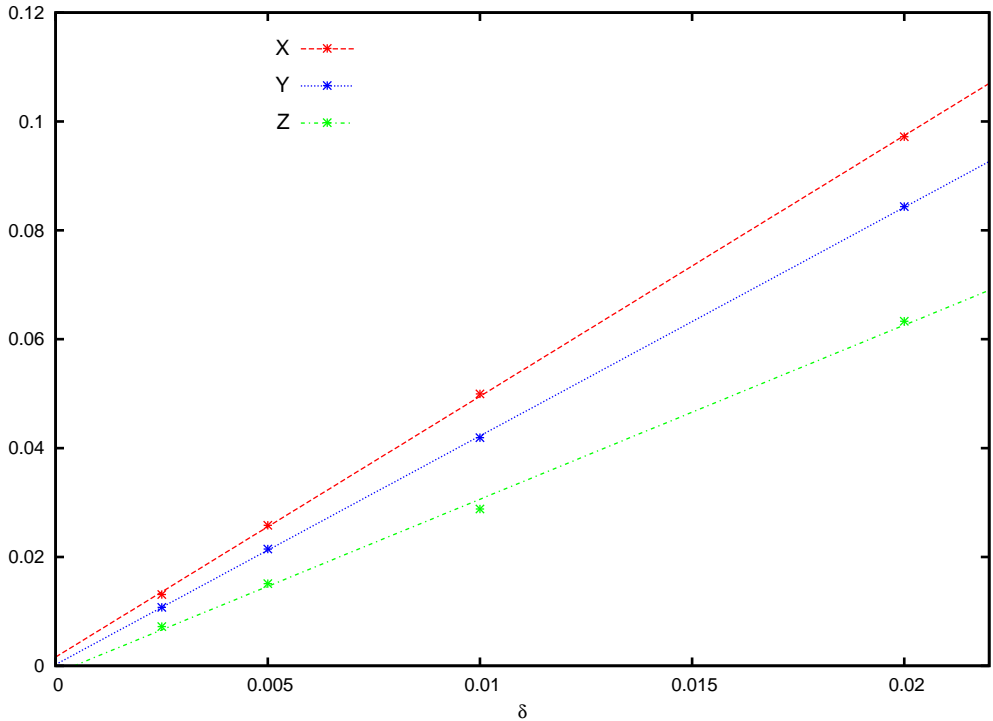}
\caption{\leftskip=25 pt \rightskip= 25 pt 
For the random variables $X,Y$ and $Z$, the plot shows the sum over the 
six subdomains of the boundary of 
the $L^1$ norms of the differences.
}
\label{fig_rt_rv}
\end{figure}

The second type of simulation in this section studies the random 
variables $X,Y,Z$ from the previous section. We divide the boundary of 
the disc into six disjoint subsets and condition on the event that the 
SKSAW exits the disc in one of these subsets. For the walks 
that meet this conditioning we rotate the lattice so that the walk 
exits at $1$. Thus for each of the three random variables we have 
six different distributions corresponding to the subset that we 
condition on. In the scaling limit we expect that they will converge
to the same limit thanks to the rotation invariance.
But before the scaling limit these six distributions are 
clearly different. 

The six subsets of the boundary of the disc are defined as follows. 
Let $\theta$ be the polar angle of the point where the walk exits the 
disc. For $\theta \in [0,\pi/4]$ the subsets are defined by 
dividing $[0,\pi/4]$ into six subintervals of equal length. 
The rotation and reflection symmetries of the square lattice
imply that any $\theta$ outside of $[0,\pi/4]$ is equivalent to 
an angle in this interval. We extend the six subsets to all of 
$[0,2\pi]$ accordingly.

For each random variable and each subset of the boundary we compute 
the $L^1$ norm of the difference between the CDF of the random variable 
conditioned on this subset and the predicted CDF. 
Then for each random variable we sum the six norms corresponding 
to the six subsets. This gives an overall measure of the deviation of the 
CDF's with the conditioning from the prediction.
Simulations were done for lattice 
spacings of $\delta=0.02, 0.01, 0.005, 0.0025$. For each $\delta$ 
we generated $10^9$ samples. So for each subset of the boundary there 
are approximately $1/6$ as many samples for the conditioned CDF. 
These sums of norms for the three random variables 
are plotted as a function of $\delta$ in figure 
\ref{fig_rt_rv}. 
Again, the lines are least squares fits to the data. 

\section{Conclusions and open questions}

We have carried out four types of simulations to test the conjecture
that the scaling limit of the smart kinetic self-avoiding 
walk on the square lattice in the full-plane is 
full-plane SLE$_6$. In one set of tests we compute the exit 
distribution for three particular domains
under the assumption that the scaling limit is rotationally invariant.
We conjecture not only that this distribution converges to 
harmonic measure, the SLE$_6$ prediction for the distribution, 
as the lattice spacing goes to zero, but also that the leading order
correction is proportional to $\delta$.
Our simulations provide excellent support
for this conjecture for both the square and hexagonal lattices.

Our second set of tests look at the full-plane 
SKSAW up to the time it exits the unit
disc centered at the origin. We assume that the scaling limit is 
rotationally invariant so that 
we can rotate the disc to put the exit point at 
$1$. We then compute the distribution of the three random variables 
$X,Y$ and $Z$. Their distributions for the analogous process using 
full-plane SLE$_6$ can be computed explicitly and we find excellent agreement 
with these predictions. Again we conjecture that the leading order
correction is proportional to $\delta$ and find strong support for 
this prediction for the square lattice.

The third set of tests is the same as the first set of tests 
except that we do not average over orientations of the lattice. 
The fourth set of tests is similar to the second set of tests.
In the second set we rotated so that the walk exits the disc at $1$. 
For the fourth set we divide the boundary of the disc into 
six subsets and condition on the event that the walk exits 
the disc in one of them, and then rotate the disc so that 
the walk exits at 1. The third and fourth sets of tests 
find strong support for the rotational invariance of the scaling
limit. 

We expect that the SKSAW is in the universality 
class of the so called $\Theta$-point for polymers. In particular, 
we expect it to be in the same universality class as the 
interacting self-avoiding walk given by \reff{isaw_ham} at its 
critical point. So if we take the $N$-step ISAW in the full-plane, 
let $N \ra \infty$ and then let $\delta \ra 0$, then we expect to 
get full-plane SLE$_6$. 
If we want to define the ISAW in a simply connected domain so that 
the scaling limit is radial or chordal SLE$_6$, then the definition
is more subtle.
The SKSAW on the hexagonal lattice helps explain what the definition 
should be. Let $D$ be a simply connected domain containing the origin.
Consider a full-plane SKSAW from $0$ to $\infty$ stopped when it 
first exits the domain. Let $\omega$ be the SKSAW up to this exit. 
It can be shown that on the hexagonal lattice the probability of $\omega$ 
essentially satisfies
\beann
P(\omega) \propto 2^{-N(\omega)}
\eeann
where $N(\omega)$ is the number of hexagons which have an edge in $\omega$. 
Thus the SKSAW up to the time it exits the domain is equivalent to 
taking all self-avoiding walks in the domain which start at the origin
and end on the boundary and weighting them according to the above equation.
The analogous construction for the ISAW is to take 
all self-avoiding walks in the domain which start at the origin
and end on the boundary and weight them by $e^{-\beta H(\omega)}$ with 
$H(\omega)$ given by \reff{isaw_ham}.
So the scaling limit of this variant of the ISAW is not radial SLE$_6$, 
but rather full-plane SLE$_6$ stopped when it exits the domain.
The energy $H(\omega)$ does not depend on the boundary.
The only effect of the 
boundary is that we only consider walks in $D$ from $0$ to the boundary.

Now fix a point $v$ on the boundary of the domain and consider the 
radial SKSAW from $0$ to $v$. 
It can be shown that on the hexagonal lattice the probability of
a walk $\omega$ in this model essentially satisfies
\beann
P_{radial}(\omega) \propto 2^{-\overline{N}(\omega)}
\eeann
where $\overline{N}(\omega)$ is the number of hexagons which have an
edge which belongs to $\omega$ 
{\it not including} the hexagons on the boundary. Let $\partial N(\omega)$
be the number of boundary hexagons that have an edge in $\omega$. So 
$\overline{N}(\omega)=N(\omega)- \partial N(\omega)$, and thus 
$P_{radial}(\omega)$ is proportional to $2^{\partial N(\omega)} P(\omega)$. 
So the radial SKSAW can be thought of as the full-plane SKSAW conditioned
to exit the domain at $v$ and weighted so that it is attracted to 
the boundary.
Thus to define a version of the ISAW that would have radial SLE$_6$
as its scaling limit, we would need to include some sort of attractive
interaction with the boundary. Presumably we would also need to 
tune the parameters in the boundary interaction to put the system 
at an appropriate critical point.

\bigskip

{\it Acknowledgments:} 
An allocation of computer time from the UA Research Computing High Performance 
Computing (HPC) and High Throughput Computing (HTC) 
at the University of Arizona is gratefully acknowledged.


\bigskip

\begin{thebibliography}{}

\bibitem{amit_et_al}
D. J. Amit, G. Parisi, L. Peliti,
Asymptotic behavior of the ``true'' self-avoiding walk.
Phys. Rev. B {\bf 27}, 1635 (1983).

\bibitem{camia_newman}
F. Camia, C. M. Newman, 
Critical percolation exploration path and SLE$_6$ : a proof of convergence.
Probab. Theory Related Fields
{\bf 139},473--519 (2007).
Archived as arXiv:math/0605035 [math.PR].

\bibitem{coniglio1987}
A. Coniglio, N. Jan, I. Majid, H. E. Stanley, 
Conformation of a polymer chain at the Theta$^\prime$ point: 
connection to the external perimeter of a percolation cluster.
Phys. Rev. B {\bf 35}, 3617 (1987).

\bibitem{duplantier_saleur}
B. Duplantier, H. Saleur, 
Exact tricritical exponents for polymers at the FTHETA point in two dimensions.
Phys. Rev. Lett. {\bf 59}, 539 (1987).


\bibitem{gherardi}
M. Gherardi, 
Theta-point polymers in the plane and Schramm-Loewner evolution.
Phys. Rev. E {\bf 88}, 032128 (2013).
Archived as arXiv:1306.4993 [cond-mat.stat-mech].

\bibitem{gunn_ortuno}
J. M. F. Gunn, M. Ortu\~no, 
Percolation and motion in a simple random environment.
J. Phys. A {\bf 18}, L1095 (1985).

\bibitem{jiang} J. Jiang,
Exploration processes and SLE$_6$. Preprint (2014). 
Archived as arXiv:1409.6834 [math.PR].

\bibitem{Kennedya} 
T.~Kennedy, 
Monte Carlo tests of SLE predictions for 2D self-avoiding walks. 
Phys. Rev. Lett. {\bf 88}, 130601 (2002).
Archived as arXiv:math/0112246v1 [math.PR]. 

\bibitem{Kennedyb} 
T.~Kennedy,	
Conformal invariance and stochastic Loewner evolution predictions 
for the 2D self-avoiding walk - Monte Carlo tests.
J. Stat. Phys. {\bf 114}, 51--78 (2004).
Archived as arXiv:math/0207231v2 [math.PR].

\bibitem{kremer1985IGSAW}
K. Kremer, J. W. Lyklema, 
Indefinitely growing self-avoiding walk.
Phys. Rev. Lett. {\bf 54}, 267 (1985).

\bibitem{lawler_book} 
G. Lawler,  {\it Conformally Invariant Processes in the Plane}.  
American Mathematical Society (2005).

\bibitem{lawler2006laplacian}
G. Lawler, 
The Laplacian-$b$ random walk and the Schramm-Loewner evolution.
Illinois J. Math. {\bf 50}, 701--746 (2006). 

\bibitem{lsw_saw} 
G.F.~Lawler, O.~Schramm, and W.~Werner, 
On the scaling limit of planar self-avoiding walk, 
{\it Fractal Geometry and Applications: a Jubilee of Benoit Mandelbrot, 
Part 2}, 339, {\em Proc. Sympos. Pure Math. 72}, 
Amer. Math. Soc., Providence, RI, 2004.
Archived as arXiv:math/0204277v2 [math.PR].

\bibitem{lerw_sle}
G. Lawler, O. Schramm, W. Werner, 
Conformal Invariance of Planar Loop-Erased Random Walks and 
Uniform Spanning Trees. 
Ann. Probab. {\bf 32}, 939--995, (2004).
Archived as arXiv:math/0112234 [math.PR].

\bibitem{lyklema1986laplacian}
J. W. Lyklema, C. Evertsz, L. Pietronero,
The Laplacian random walk.
Europhys. Lett. {\bf 2}, 77 (1986). 

\bibitem{madras_slade}
N. Madras and G. Slade, {\em The Self-Avoiding Walk}.  
Birkh\"{a}user (1996). 

\bibitem{majid1984kinetic}
I. Majid, N. Jan, A. Coniglio, H. E. Stanley,
Kinetic growth walk: A new model for linear polymers.
Phys. Rev. Lett.  {\bf 52}, 1257 (1984).

\bibitem{manna1989kinetic}
S.S. Manna, A. J. Guttmann,
Kinetic growth walks and trails on oriented square lattices: 
Hull percolation and percolation hulls.
J. Phys. A {\bf 22}, 3113 (1989).

\bibitem{smirnov}
S. Smirnov, Critical percolation in the plane: Conformal invariance, 
Cardy's formula, scaling limits.
C. R. Math. Acad. Sci. Paris
{\bf 333}, 239--244 (2001).
Archived as arXiv:0909.4499 [math.PR].

\bibitem{weinrib_trugman}
A. Weinrib, S. A. Trugman, 
A new kinetic walk and percolation perimeters.
Phys. Rev. B  {\bf 31}, 2993 (1985).

\bibitem{werner2007lectures}
W. Werner,
Lectures on two-dimensional critical percolation, 
{\it Statistical Mechanics (IAS/Park City mathemematics series v. 16)}, 
S. Sheffield, T. Spencer (eds.) (2007).
Archived as arXiv:0710.0856 [math.PR]

\bibitem{wu}
H. Wu, Conformal restriction: the radial case. Preprint (2013). 
Archived as arXiv:1304.5712 [math.PR].

\bibitem{ziff1984}
R. M. Ziff, P. T. Cummings, G. Stell,
Generation of percolation cluster perimeters by a random walk.
J. Phys. A {\bf 17},  3009  (1984).

\end{thebibliography}
\end{document}